\documentclass[11pt,oneside,english]{amsart}
\usepackage[T1]{fontenc}
\usepackage[latin9]{inputenc}
\usepackage{geometry}
\usepackage[active]{srcltx}
\usepackage{textcomp}
\usepackage{amstext}
\usepackage{amsthm}
\usepackage{amssymb}

\makeatletter
\numberwithin{equation}{section}
\numberwithin{figure}{section}
\theoremstyle{plain}
\newtheorem{thm}{\protect\theoremname}[section]
\theoremstyle{remark}
\newtheorem{rem}[thm]{\protect\remarkname}
\theoremstyle{plain}
\newtheorem{conjecture}[thm]{\protect\conjecturename}
\theoremstyle{plain}
\newtheorem{lem}[thm]{\protect\lemmaname}

\def\makebbb#1{
    \expandafter\gdef\csname#1\endcsname{
        \ensuremath{\Bbb{#1}}}
}\makebbb{R}\makebbb{N}\makebbb{Z}\makebbb{C}\makebbb{H}\makebbb{E}\makebbb{H}\makebbb{P}\makebbb{B}\makebbb{Q}\makebbb{E}\makebbb{E}
\makebbb{H}
\usepackage{babel}

\usepackage{babel}

\makeatother

\usepackage{babel}

\makeatother

\usepackage{babel}
\providecommand{\conjecturename}{Conjecture}
\providecommand{\lemmaname}{Lemma}
\providecommand{\remarkname}{Remark}
\providecommand{\theoremname}{Theorem}

\begin{document}
\title{Emergent complex geometry }
\author{Robert J. Berman}
\begin{abstract}
This is a double exposure of the probabilistic construction of Kähler-Einstein
metrics on a complex projective algebraic variety $X$ - where the
Kähler-Einstein metric emerges from a canonical random point process
on $X$ - and the variational approach to the Yau-Tian-Donaldson conjecture,
highlighting their connections. The final section is a report on joint
work in progress with Sébastien Boucksom and Mattias Jonsson on how
the non-Archimedean geometry of $X$ (with respect to the trivial
absolute value) also emerges from the probabilistic framework. 
\end{abstract}

\address{Robert J. Berman, Mathematical Sciences, Chalmers University of Technology
and the University of Gothenburg, SE-412 96 Göteborg, Sweden}
\email{robertb@chalmers.se}
\maketitle

\section{Introduction}

A recurrent theme in geometry is the quest for canonical metrics on
a given manifold $X.$ The prototypical case is when $X$ is a compact
orientable two-dimensional surface, which can be endowed with a metric
of constant scalar curvature, essentially uniquely determined by a
complex structure $J$ on $X.$ On the other hand, from a physical
point of view geometrical shapes - as we know them from everyday experience
- are, of course, not fundamental physical entities. They merely arise
as macroscopic \emph{emergent }features of ensembles of microscopic
point particles in the limit as the number $N$ of particles tends
to infinity. In mathematical terms such microscopical ensembles are\emph{
random point processes}, i.e. they are represented by a probability
measure on the configuration space of $N$ points on $X,$ or equivalently:
a symmetric probability measure $\mu^{(N)}$ on the $N-$fold product
$X^{N}.$ One is thus led to ask whether a given manifold $X$ may
be endowed with a canonical random point process - defined without
reference to any metric - from which a canonical metric $g$ emerges
as $N\rightarrow\infty?$ Here we shall focus on Kähler metrics with
constant Ricci curvature. From the physics perspective these arise
as solutions to Einstein's equations in vacuum (with Euclidean signature).
The Kähler condition means that $X$ is compatible with an integrable
complex structure $J$ on $X$ (in that parallel translation preserves
the complex structure $J).$ Such metrics - known as\emph{ Kähler-Einstein
metrics} - play a central role in current complex geometry and the
study of complex algebraic varieties, in particular in the context
of the Yau-Tian-Donaldson conjecture \cite{do2} and the Minimal Model
Program in birational algebraic geometry \cite{ko}. When a projective
algebraic variety $X$ admits a Kähler-Einstein metric it is essentially
unique, i.e. canonically attached to $X$ and can thus be leveraged
to probe $X$ using differential-geometric techniques (as, for example,
in the construction of moduli spaces \cite{s}). 

One virtue of the probabilistic approach is that it leads to essentially
explicit period type integral formulas for canonical Kähler metrics
converging towards the Kähler-Einstein metric as $N\rightarrow\infty$
(see formula \ref{eq:def of omega k KX pos intro}). These formulas
are reminiscent of the few explicit formulas for Kähler-Einstein metrics
that are available on special complex curves, involving hypergeometric
integrals (notably the modular curve, the Klein curve and Fermat curves;
see \cite[Section 2.1]{berm11}). The probabilistic approach also
generates new connections between Kähler geometry and algebraic geometry
in the context of the Yau-Tian-Donaldson conjecture on Fano varieties,
through the concept of Gibbs stability and the related stability threshold
$(\delta-$invariant) \cite{f-o,bl-j}. The present contribution to
the 2022 ICM proceedings attempts a double exposure of the probabilistic
approach in \cite{berm8,berm8 comma 5,berm6} and the variational
approach to the Yau-Tian-Donaldson conjecture in \cite{bbj}, highlighting
their connections. For more details and background the reader is referred
to the survey \cite{berm11} . See also \cite{b-c-p} for connections
between the present probabilistic approach to Kähler geometry and
quantum gravity in the context of the AdS/CFT correspondence and \cite{berm10,du}
for connections to polynomial approximation theory and pluripotential
theory in $\C^{n}.$ 

\subsection{Acknowledgments}

It is a great pleasure to thank Bo Berndtsson, Sébastien Boucksom,
Tamas Darvas, Philippe Eyssidieux, Vincent Guedj, Mattias Jonsson,
Chinh Lu, David Witt-Nyström and Ahmed Zeriahi for the stimulating
colaborations that paved the way for the work exposed here. Also thanks
to Sébastien Boucksom, Jakob Hultgren and Mingchen Xia for helpful
comments on a draft of the present manuscript.

\section{\label{sec:Emergent-K=0000E4hler-geometry}Emergent Kähler geometry}

Let $X$ be a compact complex manifold, whose dimension over $\C$
will be denoted by $n.$ The existence of a Kähler-Einstein metric
$\omega_{KE}$ on $X,$ i.e. a Kähler metric with constant Ricci curvature,
\begin{equation}
\text{Ric \ensuremath{\omega}=\ensuremath{-\beta}\ensuremath{\omega}},\label{eq:KE eq}
\end{equation}
 implies that the\emph{ canonical line bundle $K_{X}$} of $X$ (the
top exterior power of the cotangent bundle of $X)$ has a definite
sign, when $\beta\neq0$, 
\begin{equation}
\text{sign}(K_{X})=\text{sign}(\beta).\label{eq:sing of K X}
\end{equation}
We will be using the standard terminology of positivity in complex
geometry: a line bundle is said to be \emph{positive,} $L>0$ if $L$
carries some Hermitian metric with strictly positive curvature (or
equivalently, $L$ is ample in the algebro-geometric sense). The standard
additive notation for tensor products of line bundles will be adopted.
Accordingly, the dual of $L$ is expressed as $-L$ and $L$ is thus
said be negative, $L<0$ if $-L>0.$ In general, when $\beta\neq0$
the manifold $X$ is automatically a complex projective algebraic
manifold and after a rescaling of the Kähler-Einstein metric we may
as well assume that $\beta=\pm1.$ For example, in the case when $X$
is a hypersurface in $\P_{\C}^{n+1},$ cut out by a homogeneous polynomial
of degree $d,$ $K_{X}>0$ when $d>n+2,$ and $-K_{X}>0$ when $d<n+2.$ 
\begin{rem}
\label{rem:log}In the more general ``logarithmic'' setup $X$ is
replaced by a\emph{ log pair} $(X,\Delta)$ consisting of a $\Q-$divisor
$\Delta$ on a normal variety $X$ and $K_{X}$ is replaced by $K_{X}+\Delta,$
assumed to be a $\Q-$line bundle. The corresponding log Kähler-Einstein
equation \ref{eq:KE eq} is obtained by replacing $\text{Ric \ensuremath{\omega} }$with
$\text{Ric \ensuremath{\omega}\ensuremath{-[\Delta],}}$where $[\Delta]$
denotes the current of integration corresponding to $\Delta.$ For
simplicity we will stick to the case when $X$ is non-singular and
$\Delta$ is trivial (but all the results surveyed in this and the
following section generalize to the logarithmic setting, assuming
that $(X,\Delta)$ is klt \cite{berm8 comma 5,berm12,bbegz}).
\end{rem}

Coming back to the question of emergence of geometry, discussed in
the introduction, a Kähler-Einstein metric $g_{KE}$ has the crucial
property that it can be readily recovered from its volume form $dV_{KE}$,
in the case $\beta\neq0.$ Indeed, in local terms $g_{KE}$ is proportional
to the complex Hessian of the logarithm of the local density of $dV_{KE}$
(see formula\ref{eq:formula for Ricci}). Thus in order to probalistically
construct the Kähler-Einstein metric one just needs to construct a
random point process on $X$ with $N$ particles such that the empirical
measure
\begin{equation}
\delta_{N}:=\frac{1}{N}\sum_{i=1}^{N}\delta_{x_{i}},\label{eq:def of empirical measure}
\end{equation}
 viewed as a random discrete probability measure on $X,$ converges
in probability towards $dV_{KE},$ as $N\rightarrow\infty,$ 

\subsection{The case $K_{X}>0$ $(\beta=1)$}

The starting point for the probabilistic approach is the observation
that there is a canonical symmetric probability measure $\mu^{(N)}$
on the $N-$fold product $X^{N}$ of $X.$ More precisely, the integers
$N$ are taken to be of the special form
\[
N=N_{k}:=\dim_{\C}H^{0}(X,kK_{X}),
\]
 where $H^{0}(X,kK_{X})$ denotes the complex vector space of all
holomorphic section of the $k$ th tensor power of the canonical line
bundle $K_{X}\rightarrow X.$ Recall that the elements $s^{(k)}$
of $H^{0}(X,kK_{X})$ are called pluricanonical forms and may be represented
by local holomorphic functions transforming as $dz^{\otimes k},$
in terms of local holomorphic coordinates $z\in\C^{n}$ on $X.$ As
a consequence, $\left|s^{(k)}(z)\right|^{2/k}$ transforms as a local
density on $X$ and thus defines a global measure on $X.$ Replacing
$X$ with $X^{N_{k}}$ the canonical symmetric probability measure
$\mu^{(N_{k})}$ on $X^{N_{k}}$ is now defined by 
\begin{equation}
\mu^{(N_{k})}=\frac{1}{\mathcal{Z}_{N_{k}}}\left|\det S^{(k)}\right|^{2/k},\,\,\,\mathcal{Z}_{N_{k}}:=\int_{X^{N_{k}}}\left|\det S^{(k)}\right|^{2/k}\label{eq:canonical prob for K X}
\end{equation}
 where $\det S^{(k)}$ is the holomorphic section of the line bundle
$(kK_{X^{N_{k}}})\rightarrow X^{N_{k}},$ expressed as the Slater
determinant
\begin{equation}
(\det S^{(k)})(x_{1},x_{2},...,x_{N}):=\det(s_{i}^{(k)}(x_{j})),\label{eq:slater determinant}
\end{equation}
 in terms of a given basis $s_{i}^{(k)}$ in $H^{0}(X,kK_{X}).$ Under
a change of bases the section $\det S^{(k)}$ only changes by a multiplicative
complex constant (the determinant of the change of bases matrix on
$H^{0}(X,kK_{X})$).  As a consequence, $\mu^{(N_{k})}$ is independent
of the choice of bases in $H^{0}(X,kK_{X})$ and since $\det S^{(k)}$
is anti-symmetric this means that the probability measure $\mu^{(N_{k})}$
indeed defines a canonical symmetric probability measure on $X^{N_{k}}.$
Moreover, it is completely encoded by algebro-geometric data in the
following sense: realizing $X$ as projective algebraic subvariety
the section $\det S^{(k)}$ can be identified with a homogeneous polynomial,
determined by the coordinate ring of $X.$

The assumption that $K_{X}>0$ ensures that $N_{k}\rightarrow\infty$
as $k\rightarrow\infty.$ To simplify the notation we will often drop
the subindex $k$ on $N_{k}$ and consider the large $N-$limit. The
following convergence result was shown in \cite{berm8}:
\begin{thm}
\label{thm:ke conv KX pos}Let $X$ be a compact complex manifold
with positive canonical line bundle $K_{X}.$ Then the empirical measures
$\delta_{N}$ of the corresponding canonical random point processes
on $X$ (formula \ref{eq:def of empirical measure}) converge in probability,
as $N\rightarrow\infty,$ towards the normalized volume form $dV_{KE}$
of the unique Kähler-Einstein metric $\omega_{KE}$ on $X.$ 
\end{thm}

In fact, the proof  shows that the convergence holds at an exponential
rate, in the sense of large deviation theory: for any given $\epsilon>0$
there exists a positive constant $C_{\epsilon}$ such that 
\begin{equation}
\text{Prob}\left(d\left(\frac{1}{N}\sum_{i=1}^{N}\delta_{x_{i}},dV_{KE}\right)>\epsilon\right)\leq C_{\epsilon}e^{-N\epsilon},\label{eq:exp conv}
\end{equation}
 where $d$ denotes any metric on the space $\mathcal{P}(X)$ of probability
measures on $X$ compatible with the weak topology. The convergence
in probability in the previous theorem implies, in particular, that
the measures $dV_{k}$ on $X,$ defined by the expectations $\E(\delta_{N_{k}})$
of the empirical measure $\delta_{N_{k}}$ converge towards $dV_{KE}$
in the weak topology of measures on $X.$ Concretely, $dV_{k}$ is
obtained by integrating $\mu^{(N_{k})}$ over the fibers of the projection
from $X^{N_{k}}$ onto the first factor $X:$ 
\[
dV_{k}:=\int_{X^{N_{k}-1}}\mu^{(N_{k})}\rightarrow dV_{KE},\,\,\,k\rightarrow\infty
\]
For $k$ sufficiently large (ensuring that $kK_{X}$ is very ample)
the measures $dV_{k}$ are, in fact, volume forms on $X$ and induce
a sequence of canonical Kähler metrics $\omega_{k}$ on $X$ \cite[Prop 5.3]{berm8 comma 5}: 

\begin{equation}
\omega_{k}:=\frac{i}{2\pi}\partial\bar{\partial}\log dV_{k}=\frac{i}{2\pi}\partial\bar{\partial}\log\int_{X^{N_{k}-1}}\left|\det S^{(k)}\right|^{2/k}\label{eq:def of omega k KX pos intro}
\end{equation}
The convergence above also implies that the canonical Kähler metrics
$\omega_{k}$ converge, as $k\rightarrow\infty,$ towards the Kähler-Einstein
metric $\omega_{KE}$ on $X,$ in the weak topology. More generally,
as shown in \cite{berm8 comma 5}, the convergence holds on any variety
$X$ of positive Kodaira dimension (i.e. such that $N_{k}\rightarrow\infty,$
as $k\rightarrow\infty)$ if $dV_{KE}$ and $\omega_{KE}$ are replaced
by the canonical measure and current on $X,$ respectively, introduced
by Song-Tian and Tsuji in different geometric contexts \cite{berm8 comma 5}
(in the case when $X$ is singular it is assumed that $X$ is klt
and $k$ is assumed to be sufficently divisible to ensure that $kK_{X}$
is a bona fide line bundle).

\subsection{The Fano case, $K_{X}<0$ ($\beta=-1$)}

When $-K_{X}$ is positive, which means that $X$ is a Fano manifold,
any Kähler-Einstein metric on $X$ has positive Ricci curvature. However,
not all Fano manifolds $X$ carry Kähler-Einstein-metrics; according
to the Yau-Tian-Donaldson conjecture (discussed in Section \ref{sec:The-Yau-Tian-Donaldson-conjectur})
a Fano manifold admits a Kähler-Einstein-metric if and only if $X$
is K-polystable. In the probabilistic approach a new type of stability
assumption naturally appears, as next explained. First note that when
$-K_{X}>0$ the spaces $\dim H^{0}(X,kK_{X})$ are trivial for all
positive integers $k.$ On the other hand, the dimensions tend to
infinity as $k\rightarrow-\infty.$ Thus it is natural to replace
$k$ with $-k$ in the previous constructions. In particular, given
a positive integer $k,$ we set 
\[
N_{k}:=\dim H^{0}(X,-kK_{X})
\]
 and attempt to define a probability measure on $X^{N_{k}}$ as
\[
\mu^{(N_{k})}:=\frac{\left|\det S^{(k)}\right|^{-2/k}}{\mathcal{Z}_{N_{k}}},\,\,\,\,\mathcal{Z}_{N_{k}}:=\int_{X^{N_{k}}}\left|\det S^{(k)}\right|^{-2/k},
\]
 where the numerator defines a measure on the complement in $X^{N_{k}}$
of the zero-locus of $\det S^{(k)}.$ However, it may happen that
the normalizing constant $\mathcal{Z}_{N_{k}}$ diverges, since the
integrand of $\mathcal{Z}_{N_{k}}$ blows-up along the zero-locus
in $X^{N_{k}}$ of $\det S^{(k)}.$ Accordingly, a Fano manifold $X$
is called\emph{ Gibbs stable at level $k$} if $Z_{N_{k}}<\infty$
and \emph{Gibbs stable} if it is Gibbs stable at level $k$ for $k$
sufficiently large. We thus arrive at the following probabilistic
analog of the Yau-Tian-Donaldson conjecture posed in \cite{berm8 comma 5}:
\begin{conjecture}
\label{conj:Fano with triv autom intr}Let $X$ be Fano manifold.
Then 
\begin{itemize}
\item $X$ admits a unique Kähler-Einstein metric $\omega_{KE}$ if and
only if $X$ is Gibbs stable. 
\item If $X$ is Gibbs stable, the empirical measures $\delta_{N}$ of the
corresponding canonical point processes converge in probability towards
the normalized volume form of $\omega_{KE}.$ 
\end{itemize}
\end{conjecture}

It should be stressed that the Gibbs stability of $X$ implies that
the group $\text{Aut \ensuremath{(X)}}$ of automorphisms of $X$
is finite \cite[Prop 6.5]{berm8 comma 5}. Accordingly, when comparing
Conjecture \ref{conj:Fano with triv autom intr} with the Yau-Tian-Donaldson
conjecture one should view  Gibbs stability as the analog of K-stability.
There is also a natural analog of the stronger notion of uniform K-stability
\cite{bhj1,der}. To see this first note that Gibbs stability can
be given a purely algebro-geometric formulation, saying that the $\Q-$divisor
$D_{N_{k}}$ in $X^{N_{k}}$ cut out by the (multi-valued) holomorphic
section $(\det S^{(k)})^{1/k}$ of $-K_{X^{N_{k}}}$ has mild singularities
in the sense of the Minimal Model Program \cite{ko0}. More precisely,
$X$ is Gibbs stable at level $k$ iff $D_{N_{k}}$ is\emph{ klt (Kawamata
Log Terminal).} This means that the \emph{log canonical threshold
(lct) }of $D_{N_{k}}$ satisfies $\text{lct}\ensuremath{(D_{N_{k}})>1,}$
as follows directly from the standard analytic representation of the
log canonical threshold of a $\Q-$divisor as an integrability threshold
\cite{ko0}. Accordingly, $X$ is called\emph{ uniformly Gibbs stable
}if the there exists $\epsilon>0$ such that, for $k$ sufficiently
large, $\text{lct}\ensuremath{(D_{N_{k}})>1+\epsilon.}$ One is thus
led to pose the following purely algebro-geometric conjecture:
\begin{conjecture}
\label{conj:unif Gibbs iff unif K}Let $X$ be a Fano manifold. Then
$X$ is (uniformly) K-stable iff $X$ is (uniformly) Gibbs stable.
\end{conjecture}

One direction of the uniform version of the previous conjecture was
established in \cite{f-o,fu2}, using techniques from the Minimal
Model Program: 
\begin{thm}
\label{thm:(Fujita-Odaka)-Uniform-Gibbs} \cite{f-o}Uniform Gibbs
stability implies uniform K-stability
\end{thm}

Let us briefly recall the elegant argument in \cite{f-o}, which introduces
the invariant $\delta(X),$ which has come to play a key role in recent
developments around the Yau-Tian-Donaldson conjecture. First, by \cite[Thm 2.5]{f-o},
\begin{equation}
\text{lct }(D_{N_{k}})\leq\delta_{k}(X):=\inf_{\Delta_{k}}\text{lct }(\Delta_{k}),\label{eq:def of delta k}
\end{equation}
 where the inf is taken over all anti-canonical $\Q-$divisors $\Delta_{k}$
on $X$ of $k-$basis type, i.e. $\Delta_{k}$ is the normalized sum
of the $N_{k}$ zero-divisors on $X$ defined by the members of a
given basis in $H^{0}(X,-kK_{X})$. Finally, by \cite[Thm 0.3]{f-o},
if the invariant $\delta(X)$ defined as 
\begin{equation}
\delta(X):=\limsup_{k\rightarrow\infty}\delta_{k}(X)\label{eq:def of delta as limit of delta k}
\end{equation}
 satisfies $\delta(X)>1,$ then $X$ is uniformly K-stable \cite{fu2}
and thus admits a unique Kähler-Einstein metric by the solution of
the (uniform) Yau-Tian-Donadlson conjecture recalled in Section \ref{subsec:The-variational-approach}.
In particular, this means that uniform Gibbs stability implies the
existence of a Kähler-Einstein metrics (in line with Conjecture \ref{conj:Fano with triv autom intr}).
However, the converse implication, that we shall come back to in Section
\ref{sec:A-Non-Archimedean-approach}, is still open. Anyhow, even
if confirmed, it is a separate analytic problem to prove the convergence
in Conjecture \ref{conj:Fano with triv autom intr}. ``Tropicalized''
analogs of Conjecture \ref{conj:Fano with triv autom intr} are established
on toric varieties in \cite{b-oe} and on tori in \cite{hu}.

In \cite{berm11} a variational approach to the convergence problem
was introduced, further developed in \cite{berm12}, where the convergence
was settled on log Fano curves. In the general case the approach yields,
in particular, the following conditional convergence result:
\begin{thm}
\label{thm:conv beta neg}\cite{berm11,berm12} Let $X$ be a Fano
manifold and assume that $X$ admits a Kähler-Einstein metric $\omega_{KE}.$
Take the basis $s_{i}^{(k)}$ in formula \ref{eq:slater determinant}
to be orthonormal wrt the Hermitian metric on $H^{0}(X,-kK_{X})$
induced by $\omega_{KE}$ and assume that 
\begin{equation}
\lim_{N\rightarrow\infty}\frac{1}{N}\log\mathcal{Z}_{N}=0.\label{eq:log Z N conv zero}
\end{equation}
Then $\text{Aut \ensuremath{(X)}}$ is finite and the empirical measure
$\delta_{N}$ converge in probability towards the normalized volume
form $dV_{KE}$ of the unique Kähler-Einstein-metric $\omega_{KE}$
on $X.$ 
\end{thm}

In \cite{berm11} two different types of hypotheses were put forth
ensuring that the convergence \ref{eq:log Z N conv zero} holds, one
of which  will be recalled in Section \ref{subsec:The-case beta negative proof}.
The other one assumes, in particular, that the partition function
$\mathcal{Z}_{N}(\beta),$ discussed in the following section, is
zero-free in some $N-$independent neighborhood $\Omega$ of $]-1,0]$
in $\C$ (when $\mathcal{Z}_{N}(\beta)$ is analytically continued
to a holomorphic function on $\Omega).$ This allows one to ``analytically
continue'' the convergence when $\beta>0$ to $\beta<0.$ This is
discussed in detail in \cite{berm12}, where some intriguing connections
between this zero-free hypothesis and the zero-free property of the
local L-functions appearing in the Langlands program are also pointed
out. 

\subsection{\label{subsec:The-statistical-mechanical}The statistical mechanical
formalism and outlines of the proofs}

Theorem \ref{thm:ke conv KX pos} (or more precisely, the exponential
convergence in formula \ref{eq:exp conv}) is deduced from a Large
Deviation Principle (LDP), which may be symbolically expressed as
\begin{equation}
\text{Prob }\left(\frac{1}{N}\sum_{i=1}^{N}\delta_{x_{i}}\in B_{\epsilon}(\mu)\right)\sim e^{-NR(\mu)},\,\,\,N\rightarrow\infty,\,\,\epsilon\rightarrow0\label{eq:LDP}
\end{equation}
 where $B_{\epsilon}(\mu)$ denotes the ball of radius $\epsilon$
centered at a given $\mu$ in the space $\mathcal{P}(X)$ of all probability
measure on $X,$ endowed with a metric $d$ compatible with the weak
topology. In probabilistic terminology the functional $R(\mu)$ is
called the \emph{rate functional}. By general principles, any rate
functional of a LDP is lower-semicontinuous and its infimum vanishes.
In the present setup the volume form $dV_{KE}$ of the Kähler-Einstein
metric is the unique minimizer of $R(\mu),$ which yields the exponential
convergence in formula \ref{eq:exp conv}. 

As next explained, the proof of the LDP is inspired by statistical
mechanics. Fix a Kähler metric on $X$. It induces a a volume form
$dV$ on $X$ and a Hermitian metric $\left\Vert \cdot\right\Vert $
on $K_{X}.$ The canonical probability measure \ref{eq:canonical prob for K X}
may then be decomposed as 

\[
\mu^{(N)}=\frac{1}{\mathcal{Z}_{N_{k}}}\left\Vert \det S^{(k)}\right\Vert ^{2/k}dV^{\otimes N},
\]
where the basis $s_{i}^{(k)}$ in formula \ref{eq:slater determinant}
is taken to be orthonormal wrt the Hermitian metric on $H^{0}(X,kK_{X})$
induced by $dV$ and $\left\Vert \cdot\right\Vert .$ Introducing
the \emph{energy per particle }
\begin{equation}
E^{(N)}(x_{1},...,x_{N}):=-\frac{1}{kN}\log\left\Vert \det S^{(k)}(x_{1},...,x_{N_{k}})\right\Vert ^{2}\label{eq:def of E N intro}
\end{equation}
we can thus express $\mu^{(N)}$ as the following \emph{Gibbs measure,}
at\emph{ inverse temperature $\beta=1:$}
\begin{equation}
\mu_{\beta}^{(N)}=\frac{e^{-\beta NE^{(N)}}}{\mathcal{Z}_{N}(\beta)}dV^{\otimes N},\,\,\,\mathcal{Z}_{N}(\beta):=\int_{:X^{N}}e^{-\beta NE^{(N)}}dV^{\otimes N}.\label{eq:def of Gibbs measure intro}
\end{equation}
 In statistical mechanical terms the Gibbs measures represents the
microscopic thermal equilibrium state of $N$ interacting identical
particles on $X.$ The normalizing constant $\mathcal{Z}_{N}(\beta)$
is called the \emph{partition function. }

The starting point of the proof of the LDP \ref{eq:LDP} is a classical
result of Sanov in probability, going back to Boltzmann, saying that
in the ``non-interacting case'' $\beta=0$ (where the positions
$x_{i}$ define independent random variables on $X)$ the LDP holds
with rate functional given by the entropy Ent$(\mu)$ of $\mu$ relative
to $dV,$ i.e. the functional on $\mathcal{P}(X)$ defined by
\[
\text{Ent}(\mu):=\int_{X}\log(\frac{\mu}{dV})\mu
\]
if $\mu$ is absolutely continuous wrt $dV$ and otherwise $\text{Ent}(\mu):=+\infty$\footnote{In the physics literature the opposite sign convention for $\text{Ent}(\mu)$
is used. } The strategy to handle the ``interacting case'' $\beta\neq0$ is
to first show that there exists a functional $E(\mu)$ on $\mathcal{P}(X)$
such that the energy per particle, $E^{(N)}(x_{1},...x_{N}),$ may
be approximated as
\begin{equation}
E^{(N)}(x_{1},...x_{N})\rightarrow E(\mu),\label{eq:Gamma conv}
\end{equation}
when $\frac{1}{N}\sum_{i=1}^{N}\delta_{x_{i}}\rightarrow\mu,$ in
an appropriate sense, as $N\rightarrow\infty.$ Formally combining
this result with Sanov's LDP suggests that for any $\beta>0$ the
corresponding rate functional is given by 
\begin{equation}
R_{\beta}(\mu)=F_{\beta}(\mu)-\inf_{\mathcal{P}(X)}F_{\beta},\,\,\,F_{\beta}(\mu)=\beta E(\mu)+\text{Ent}(\mu)\in]0,\infty],\label{eq:rate funct and free en}
\end{equation}

In thermodynamical terms the functional $F_{\beta}(\mu)$ is the \emph{free
energy,} at inverse temperature $\beta$ (strictly speaking it is
$\beta^{-1}F_{\beta}$ which is the free energy, i.e. the energy that
is free to do work once the disordered thermal energy has been subtracted).
In the present setting the role of the ``macroscopic'' energy $E(\mu)$
is played by the \emph{pluricomplex energy }of the measure $\mu$
(introduced in \cite{bbgz} and discussed in Section \ref{sec:The-thermodynamical-formalism}).
Briefly, it is first shown in \cite{berm8} that the convergence \ref{eq:Gamma conv}
holds in the sense of \emph{Gamma-convergence. }This means that 
\begin{equation}
\frac{1}{N_{j}}\sum_{i=1}^{N_{j}}\delta_{x_{i}}\rightarrow\mu\implies\liminf_{N_{j}\rightarrow\infty}E^{(N_{j})}(x_{1},...x_{N_{j}})\geq E(\mu)\label{eq:limif for E N}
\end{equation}
 and for any $\mu$ there exists some sequence of configurations in
$X^{N}$ saturating the previous inequality. The Gamma-convergence
is deduced from the convergence of weighted transfinite diameters
established in \cite{b-b} using a duality argument (where $E(\mu)$
arises as a Legendre-Fenchel transform; compare formula \ref{eq:Leg tr of E}).
The combination with Sanov's theorem is then made rigorous using an
effective submean inequality on small balls in the Riemannian orbifold
$X^{N}/S_{N},$ established using geometric analysis. 

The free energy functional $F_{\beta}$ has a unique minimizer $\mu_{\beta}$
in $\mathcal{P}(X)$ for any $\beta>0$ (as discussed in Section \ref{subsec:Back-to-the}).
As a consequence, the empirical measure $\delta_{N}$ converges in
probability towards $\mu_{\beta},$ as $N\rightarrow\infty.$ The
LDP proved in \cite{berm8} also implies that for $\beta>0$
\begin{equation}
\lim_{N\rightarrow\infty}-\frac{1}{N}\log\mathcal{Z}_{N}(\beta)=\inf_{\mathcal{P}(X)}F_{\beta}\label{eq:conv of part funct}
\end{equation}
Incidentally, the free energy functional $F_{\beta}$ on $\mathcal{P}(X)$
may be identified with the (twisted) Mabuchi functional in Kähler
geometry, as explained in Section \ref{subsec:The-Mabuchi-and}. 

\subsubsection{\label{subsec:The-case beta negative proof}The case $\beta<0$}

The Gibbs measure $\mu_{\beta}^{(N)}$ can, alternatively, be viewed
as a Gibbs measure at \emph{unit }temperature, if $E^{(N)}$ is replaced
with with the rescaled energy $\beta E^{(N)}$ (thus treating $\beta$
as a coupling constant). For $\beta>0$ this energy is is \emph{repulsive,}
since it tends to $\infty$ as any two particle positions merge (due
to the vanishing of the determinant $\det S^{(k)}(x_{1},...,x_{N_{k}})$).
However, when $\beta$ changes sign the rescaled energy $\beta E^{(N)}$
becomes\emph{ attractiv}e; it tends to $-\infty$ as any two points
merge, which leads to subtle concentration phenomena and various new
technical difficulties. For example, one reason that the proof of
the LDP does not generalize to $\beta<0$ is that the Gamma-convergence
in formula \ref{eq:Gamma conv} is not preserved when $E^{(N)}$ is
replaced by $-E^{(N)}.$ In order to bypass this difficulty a variational
approach was introduced in \cite{berm11}. The starting point is the
classical Gibbs variational principle, which yields
\begin{equation}
-\frac{1}{N}\log\mathcal{Z}_{N}(\beta)=\inf_{\mathcal{P}(X^{N})}F_{\beta}^{(N)},\,\,\,\,\,\,F_{\beta}^{(N)}(\cdot):=\beta\left\langle E^{(N)},\cdot\right\rangle +N^{-1}\text{Ent}(\cdot),\label{eq:Gibbs var pr}
\end{equation}
 where the functional $F_{\beta}^{(N)}$ on $\mathcal{P}(X^{N})$
is called the \emph{$N-$particle mean free energy} and $\text{Ent}(\cdot)$
denote the entropy relative to $dV^{\otimes N}.$ When its infimum
is finite it is uniquely attained at the corresponding Gibbs measure
$\mu_{\beta}^{(N)}.$ In \cite{berm11,berm12} this variational formulation
is leveraged to show that, if $X$ admits a Kähler-Einstein-metric
$dV_{KE},$ then $\delta_{N}$ converges in probability towards $dV_{KE},$
under the assumption that the convergence of the partition functions
\ref{eq:conv of part funct} holds at $\beta=-1.$ In particular,
when the fixed metric on $X$ is taken to be a Kähler-Einstein metric
this proves Theorem \ref{thm:conv beta neg}, since $F_{-1}(dV_{KE})=0.$
Moreover,  the convergence \ref{eq:conv of part funct} of the partition
functions at $\beta=-1$ is shown to be implied by the following 
\begin{equation}
\text{Hypothesis:\,\,\,}\lim_{N_{j}\rightarrow\infty}(\delta_{N_{j}})_{*}\mu_{-1}^{(N_{j})}=\Gamma\in\mathcal{P}\left(\mathcal{P}(X)\right)\implies\limsup_{N_{j}\rightarrow\infty}\left\langle E^{(N_{j})},\mu_{-1}^{(N_{j})}\right\rangle \leq\left\langle E,\Gamma\right\rangle \label{eq:hyp for mean en}
\end{equation}
 where $(\delta_{N})_{*}\mu_{-1}^{(N)}$ is  the probability measure
on the infinite dimensional $\mathcal{P}(X),$ defined as the push-forward
of the canonical probability measure $\mu_{-1}^{(N)}$ on $X^{N}$
to $\mathcal{P}(X)$ under the map $\delta_{N}$ (the reversed inequality
holds for \emph{any} sequence $\mu_{N}$ in $\mathcal{P}(X^{N}),$
as follows from the inequality \ref{eq:limif for E N}). If the hypothesis
holds, then it follows that $\Gamma$ is the Dirac mass at $dV_{KE},$which
is equivalent to the convergence in Theorem \ref{thm:conv beta neg}.
In fact, as shown in \cite{berm12}, the previous hypothesis is ``almost''
equivalent to the convergence in Conjecture \ref{conj:Fano with triv autom intr}. 

Finally, we note that the conjectural extension of the formula \ref{eq:conv of part funct}
to any $\beta<0$ also suggests the following conjecture posed in
\cite{berm8} (the definition of the log canonical threshold $\text{lct }(D_{N})$
was discussed after Conjecture \ref{conj:Fano with triv autom intr}):
\begin{conjecture}
\label{conj:lim lct is Gamma}For any Fano manifold $X$ 
\begin{equation}
\lim_{N\rightarrow\infty}\text{lct }(D_{N})=\Gamma(X),\,\,\,\,\Gamma(X):=\sup_{\beta<0}\left\{ -\beta:\,\,\,\inf_{\mathcal{P}(X)}F_{\beta}>-\infty\right\} .\label{eq:def of Gamma X}
\end{equation}
\end{conjecture}

\section{\label{sec:The-thermodynamical-formalism}The thermodynamical formalism
and pluripotential theory}

The pluricomplex energy $E(\mu),$ appearing as the ``energy part''
of the free energy functional $F_{\beta}(\mu)$ in formula \ref{eq:rate funct and free en},
may be defined as the greatest lower semicontinuous extension to the
space $\mathcal{P}(X)$ of the functional whose first variation on
the subspace of volume forms is given by 
\begin{equation}
dE(\mu)=-u_{\mu},\label{eq:first var of E}
\end{equation}
 with $u_{\mu}\in C^{\infty}(X)$ denoting the solution to the following
complex Monge-Ampère equation (known as the \emph{Calabi-Yau equation})\emph{
}
\begin{equation}
MA(u)=\mu,\label{eq:cy eq}
\end{equation}
 expressed in terms of the complex Monge-Ampère measure $MA(u),$
whose definition we next recall. 

\subsection{Kähler geometry recap}

Assume given a line bundle $L$ endowed with a Hermitian metric $\left\Vert \cdot\right\Vert $
(in the present setup $L=\pm K_{X}$ and $\left\Vert \cdot\right\Vert $
is the metric on $L$ induced by a fixed Kähler metric on $X).$ Then
any smooth function $u$ on $X$ induces a metric $\left\Vert \cdot\right\Vert e^{-u/2}$
on $L,$ whose curvature form, multiplied by $i/2\pi,$ will be denoted
by $\omega_{u};$ it is a real closed two-form on $X,$ representing
the first Chern class $c_{1}(L)\in H^{2}(X,\Z)$ of $L.$ Concretely, 

\begin{equation}
\omega_{u}=\omega_{0}+\frac{i}{2\pi}\partial\bar{\partial}u_{\beta},\,\,\,\,\,\,\,\,\partial\bar{\partial}u:=\sum_{i,j\leq n}\frac{\partial^{2}u}{\partial z_{i}\partial\bar{z}_{i}}dz_{i}\wedge d\bar{z}_{j},\label{eq:def of omega u}
\end{equation}
 in terms of local holomorphic coordinates, where $\omega_{0}$ is
the normalized curvature form of the fixed metric $\left\Vert \cdot\right\Vert $
on $L.$ The complex Monge-Ampère measure $MA(u)$ is the normalized
volume form on $X$ defined by

\[
MA(u):=\omega_{u}^{n}/V,\,\,\,V:=\int_{X}\omega_{u}^{n}=\int_{X}\omega_{0}^{n}.
\]
 By the Calabi-Yau theorem there exists a smooth solution $u_{\mu}$
to the Calabi-Yau equation \ref{eq:cy eq}, uniquely determined up
to an additive constant. It has the property that $\omega_{u_{\mu}}$
is a\emph{ Kähler form. }Recall that a $J-$invariant closed real
form $\omega$ on $X$ is said to be Kähler if $\omega>0$ in the
sense that the corresponding symmetric two-tensor
\[
g:=\omega(\cdot,J\cdot)
\]
 is positive definite, i.e. defines a Riemannian metric (where $J$
denotes the complex structure on $X).$ In practice, one then identifies
the Kähler form $\omega$ with the corresponding \emph{Kähler metric}
$g.$ Likewise, the Ricci curvature of a Kähler metric $\omega$ may
be identified with the two-form
\begin{equation}
\text{Ric \ensuremath{\omega}}=-\frac{i}{2\pi}\partial\bar{\partial}\log dV.\label{eq:formula for Ricci}
\end{equation}
 where $dV$ denotes the volume form of $\omega.$ In other words,
$\text{Ric \ensuremath{\omega}}$ is the curvature of the metric on
$-K_{X}$ induced by $\omega.$ If the Kähler form $\omega$ is of
the form $\omega_{u}$ (as in formula \ref{eq:def of omega u}), then
$u$ is said to be a \emph{Kähler potential} for $\omega$ (relative
to $\omega_{0}).$ We will denote by $\mathcal{H}(X,\omega_{0})$
the space of all Kähler potentials, relative to $\omega_{0}$ and
by $\mathcal{H}(X,\omega_{0})_{0}$ the subspace of all sup-normalized
$u,$ $\sup_{X}u=0.$ The map 
\[
u\mapsto\omega_{u},\,\,\,\mathcal{H}(X,\omega_{0})_{0}\hookrightarrow c_{1}(L)
\]
 yields a one-to-one correspondence between $\mathcal{H}(X,\omega_{0})_{0}$
and the space of all Kähler forms in the first Chern class $c_{1}(L)$
of $L.$ Similarly, the Calabi-Yau theorem yields the ``Calabi-Yau
correspondence'' 
\begin{equation}
u\mapsto MA(u),\,\,\,\mathcal{H}(X,\omega_{0})_{0}\hookrightarrow\mathcal{P}(X)\label{eq:u maps to MA}
\end{equation}
 between $\mathcal{H}(X,\omega_{0})_{0}$ and the space of all volume
forms in $\mathcal{P}(X),$ where $u$ corresponds to the normalized
volume form of the Kähler metric $\omega_{u}.$ The one-form on $\mathcal{H}(X,\omega_{0})$
induced by $MA$ is exact, i.e. there exists a functional $\mathcal{E}$
on $\mathcal{H}(X,\omega_{0})$ such that
\[
d\mathcal{E}=MA,\,\,\,\,\,\,\text{i.e. }\frac{d\mathcal{E}(u+t\dot{u})}{dt}_{|t=0}=\left\langle MA(u),\dot{u}\right\rangle .
\]
(this functional is often denoted by $E$ in the literature \cite{begz},
but here we shall reserve capital letters for functionals defined
on $\mathcal{P}(X)$). The functional $\mathcal{E}(u)$ is uniquely
determined up to an additive a constant and may be explicitly defined
by
\begin{equation}
\mathcal{E}(u):=\frac{1}{V(n+1)}\sum_{j=0}^{n}\int_{X}u\omega_{u}^{j}\wedge\omega_{0}^{n-j}\label{eq:formula for e u}
\end{equation}

\subsection{\label{subsec:Pluripotential-theory-recap}Pluripotential theory
recap}

The analysis of the minimizers of $F_{\beta}$ involves some pluripotential
theory that we briefly recall. The space $PSH(X,\omega)$ of all $\omega_{0}-$psh
functions on $X$ may be defined as the closure of $\mathcal{H}(X,\omega_{0})$
in $L^{1}(X)$ (more precisely, any $u\in PSH(X,\omega)$ is the decreasing
limit of elements $u_{k}\in\mathcal{H}(X,\omega_{0})).$ The corresponding
sup-normalized subspace $PSH(X,\omega_{0})_{0}$ is compact in $L^{1}(X,\omega_{0}).$
By \cite{bbgz} the ``Calabi-Yau correspondence`` \ref{eq:u maps to MA}
extends to a correspondence between the subspace of probability measures
$\mu$ with finite energy and a subspace of $PSH(X,\omega_{0})$ denoted
by $\mathcal{E}^{1}(X,\omega_{0}):$ 

\begin{equation}
MA:\,\,\,\mathcal{E}^{1}(X,\omega_{0})_{0}\longleftrightarrow\left\{ \mu\in\mathcal{P}(X):\,\,E(\mu)<\infty\right\} \label{eq:MA correspondence for energy}
\end{equation}
 where $MA(u)$ is defined on $\mathcal{E}^{1}(X,\omega_{0})$ using
the notion of non-pluripolar products introduced in \cite{begz}.
The space $\mathcal{E}^{1}(X,\omega_{0})$ was originally introduced
in \cite{g-z}, but, as shown in \cite{bbgz}, it may also be defined
as the space of all $u\in PSH(X,\omega_{0})$ such that $\mathcal{E}(u)>-\infty,$
where $\mathcal{E}$ denotes the smallest upper semi-continuous extension
of $\mathcal{E}$ to $PSH(X,\omega_{0}).$

\subsection{\label{subsec:Back-to-the}Back to the free energy functional $F_{\beta}$}

The free energy functional $F_{\beta},$ defined in formula \ref{eq:rate funct and free en},
$F_{\beta}=\beta E+\text{Ent,}$ is lsc and convex on $\mathcal{P}(X)$
when $\beta>0$ (since both terms are). In the case when $\beta<0$
we define $F_{\beta}(\mu)$ by the same expression when $E_{\omega_{0}}(\mu)<\infty$
and otherwise we set $F_{\beta}(\mu)=\infty.$ The definition is made
so that we still have $F_{\mu}(\mu)\in]-\infty,\infty]$ with $F_{\mu}(\mu)<\infty$
iff both $E(\mu)<\infty$ and Ent$(\mu)<\infty.$ 

The following lemma follows readily from the first variation \ref{eq:first var of E}
and the formula \ref{eq:formula for Ricci} for Ricci curvature of
a Kähler metric.
\begin{lem}
\label{lem:critical pt}A volume form $\mu$ on $X$ is a critical
point of the functional $F_{\beta}$ on $\mathcal{P}(X)$ iff the
function 
\[
u_{\beta}:=\frac{1}{\beta}\log\frac{\mu}{dV}
\]
 solves the complex Monge-Ampère equation 
\begin{equation}
MA(u)=e^{\beta u}dV\label{eq:MA eq for beta in section thermo}
\end{equation}
iff $\omega_{\beta}:=\omega_{u_{\beta}}$ is a Kähler form solving
the twisted Kähler-Einstein equation 
\begin{equation}
\mbox{\ensuremath{\mbox{Ric}}\ensuremath{\omega}}+\beta\omega=\theta,\,\,\,\theta:=(\beta\ensuremath{\mp1)\omega_{0}}\label{eq:tw ke eq in text}
\end{equation}
\end{lem}

In the Fano case the previous equation coincides with Aubin's continuity
equation with ``time-parameter'' $t:=-\beta.$ When $\beta>0$ it
follows directly from the lower semicontinuity of $F_{\beta}$ on
the compact space $\mathcal{P}(X)$ that $F_{\beta}$ admits a minimizer.
\begin{thm}
\label{thm:reg and exist of min}\cite{berm6} 
\begin{itemize}
\item (regularity) Any minimizer $\mu_{\beta}$ of the functional $F_{\beta}$
on $\mathcal{P}(X)$ is a volume form and thus of the form in Lemma
\ref{lem:critical pt}
\item (existence) If $F_{\beta_{0}}$ is bounded from below for some $\beta_{0}<0,$
then for any $\beta>\beta_{0}$ the functional $F_{\beta}$ on $\mathcal{P}(X)$
admits a minimizer. In other words, if $F_{\beta}$ is \emph{coercive}
(wrt $E)$ in the sense that there exists $\epsilon>0$ and $C>0$
such that 
\begin{equation}
F_{\beta}\geq\epsilon E+C,\label{eq:coerciv of F beta}
\end{equation}
 then $F_{\beta}$ admits a minimizer.
\end{itemize}
\end{thm}

Moreover, by the Bando-Mabuchi theorem, if $\beta>-1$ the minimizer
is uniquely determined and if $\beta=-1$ it is uniquely determined
iff the automorphism group $\text{Aut \ensuremath{(X)} of \ensuremath{X}}$
is finite (see \cite{ber-bernd} for generalizations). The proof of
the previous theorem employs a duality argument, which fits naturally
into the thermodynamical formalism, when combined with pluripotential
theory and the variational approach to complex Monge-Ampère equation
developed in \cite{bbgz}. The strategy is to show that any minimizer
satisfies the Monge-Ampère equation \ref{eq:MA eq for beta in section thermo}
in the weak sense of pluripotential theory, so that the regularity
theory for Monge-Ampère equations (going back to Aubin and Yau), can
be invoked. In the case when $\beta>0$ the proof of Theorem \ref{thm:reg and exist of min}
follows from the strict convexity of $F_{\beta},$ resulting from
the convexity of $E(\mu)$ and the strict convexity of $\text{Ent}(\mu)$
on $\mathcal{P}(X),$ combined with the Aubin-Yau theorem \cite{au,y}
(showing that there exists a unique smooth solution to the equation
\ref{eq:MA eq for beta in section thermo}). The proof in the case
when $\beta<0$ exploits the \emph{Legendre-Fenchel transform.} Recall
that, in general, this transform yields a correspondence between lsc
convex functions on a locally convex topological vector space $V$
and its dual $V^{*}.$ In order to facilitate the comparison to the
standard functionals in Kähler geometry (discussed in the following
section) it will, however, be convenient to use a slightly non-standard
sign convention where a lsc convex function $f$ on $V$ corresponds
to the usc concave function $f^{*}$ on $V^{*}$ defined by

\begin{equation}
f^{*}(w):=\inf_{v\in V}\left(\left\langle v,w\right\rangle +f(v)\right).\label{eq:def of Leg-F trans}
\end{equation}
 Conversely, if $\Lambda$ is a functional on $V^{*}$ we define $\Lambda^{*}(v)$
as the lsc convex function
\[
\Lambda^{*}(v)=\sup_{w\in V^{*}}\left(-\left\langle v,w\right\rangle +f(w)\right).
\]
 We take $V$ to be the space of all signed measures $\mu$ on $X,$
so that $V^{*}=C^{0}(X).$ We can then view $E$ and $\text{Ent}$
as convex lsc functions on $V,$ which, by definition, are equal to
$\infty$ on the complement of $\mathcal{P}(X)$ in $V.$ Under the
Legendre-Fenchel transform these correspond to the usc convex functions
$E^{*}$ and $\text{Ent}^{*},$ respectively, on $C^{0}(X),$ which
turn out to be Gateaux differentiable. Indeed, by a classical result
(which follows from Jensen's inequality)
\[
\text{Ent}^{*}(u)=-\log\int e^{-u}dV
\]
Moreover, as shown in \cite{b-b,bbgz} the functional $E^{*}$ on
$C^{0}(X)$ is Gateaux differentiable and 
\begin{equation}
E^{*}(u)=\mathcal{E}(u),\,\,\,dE_{|u}^{*}=MA(u),\,\,\,\text{for \ensuremath{u\in\mathcal{H}(X,\omega_{0}).}}\label{eq:Leg tr of E}
\end{equation}

Now consider, for simplicity, the case $\beta=-1$ (the general case
is obtained by a simple scaling). It follows directly from the fact
that the Legendre-Fenchel transform is increasing and involutive that
\begin{equation}
\inf_{\mathcal{P}(X)}F_{-1}:=\inf_{\mathcal{P}(X)}\left(-E+\text{Ent}\right)=\inf_{C^{0}(X)}\left(-E^{*}+\text{Ent}^{*}\right)\label{eq:inf F is inf dual}
\end{equation}
Moreover, it readily from the definitions that
\[
F_{-1}\left(MA(u)\right)=\left(-E+\text{Ent}\right)(\,dE_{|u}^{*})\geq\left(-E^{*}+\text{Ent}^{*}\right)(u).
\]
Hence, if $\mu$ minimizes $F_{-1}$ and we express $\mu=MA(u_{\mu}),$
then $u_{\mu}$ minimizes the functional $-E^{*}+\text{Ent}^{*}$
on $C^{0}(X).$ However, in the present setup $u_{\mu}$ is not, a
priori, in $C^{0}(X),$ but only in $\mathcal{E}^{1}(X,\omega_{0}).$
This problem is circumvented using a simple approximation argument
to deduce that $u_{\mu}$ minimizes the extension of the functional
$\left(-E^{*}+\text{Ent}^{*}\right)$ to $\mathcal{E}^{1}(X,\omega_{0}).$
Finally, by the Gateaux differentiability of the functional $-E^{*}+\text{Ent}^{*}$
on $C^{0}(X)$ (or more precisely, on $\{u\}+C^{0}(X)$ for any given
$u\in\mathcal{E}^{1}(X,\omega_{0})$) it then follows that $u_{\mu}$
is a critical point of the functional $-E^{*}+\text{Ent}^{*}.$ Thus,
after perhaps adding a constant to $u_{\mu},$ it satisfies the complex
Monge-Ampère equation \ref{eq:MA eq for beta in section thermo} in
the weak sense of pluripotential theory.

The proof of the first point in Theorem \ref{thm:reg and exist of min}
can now be concluded by invoking the regularity results for pluripotential
solutions to Monge-Ampère equations (which, by \cite[Appendix B]{bbegz},
hold in the general setup of log Fano varieties). As for the second
point it is shown in \cite{berm6} by proving that any minimizing
sequence $\mu_{j}$ in $\mathcal{P}(X)$ (i.e. a sequence $\mu_{j}$
such that $F_{\beta}(\mu_{j})$ converges to the infimum of $F_{\beta})$
converges (after perhaps passing to a subsequence) to a minimizer
of $F_{\beta}.$ This is shown using a duality argument, as above.
Alternatively, as shown in \cite{bbegz} in a more general singular
context (including singular log Fano varieties), the existence of
a minimizer for $F_{\beta}(\mu)$ follows from the following result
in \cite{bbegz}:
\begin{thm}
\label{thm:(energy/entrop-compactness)}(energy/entropy compactness).
The functional $E(\mu)$ is continuous on any sublevel set $\{\text{Ent}\leq C\}\subset\mathcal{P}(X).$
As a consequence, if $F_{\beta}$ is coercive on $\mathcal{P}(X),$
then it is lower semi-continuous and thus admits a minimizer.
\end{thm}

This result has come to play a prominent role in recent developments
in Kähler geometry, as discussed in Section \ref{subsec:The-uniform-YTD}.

\subsection{\label{subsec:The-Mabuchi-and}The Mabuchi and Ding functionals }

Under the ``Calabi-Yau correspondence'' \ref{eq:u maps to MA} the
free energy functional $F_{\beta}$ on $\mathcal{P}(X)$ corresponds
to a functional $\mathcal{M_{\beta}}(u)$ on $\mathcal{E}^{1}(X,\omega_{0})$
defined by
\begin{equation}
\mathcal{M_{\beta}}(u):=F_{\beta}\left(MA(u)\right)\label{eq:M beta as F}
\end{equation}
Moreover, the functional $E(\mu)$ on $\mathcal{P}(X)$ corresponds
to the functional $E(MA(u))$ on $PSH(X,\omega_{0})$ which induces
an exhaustion function on $\mathcal{E}^{1}(X,\omega_{0})_{0},$ comparable
to $-\mathcal{E}(u),$ defining a notion of coercivity on $\mathcal{E}^{1}(X,\omega_{0})$
(in terms of the standard functionals $I$ and $J$ in Kähler geometry
$E(MA(u))=(I-J)(u)$). 

As is turns out, when restricted to $\mathcal{H}(X,\omega_{0})$ the
functional $\mathcal{M_{\beta}}(u)$ coincides with the \emph{(twisted)
Mabuchi functional.} The Mabuchi functional $\mathcal{M}$ associated
to a general polarized manifold $(X,L)$ was originally defined (up
to normalization) by the property that its first variation is proportional
to the scalar curvature of the Kähler metric $\omega_{u}$ minus the
average scalar curvature \cite{ma}. An ``energy+entropy'' formula
for $\mathcal{M},$ similar to formula \ref{eq:M beta as F}, holds
for a general polarized manifold, as first discovered in \cite{ti0}
and \cite{ch0}. Likewise, the functional on $\mathcal{E}^{1}(X,\omega_{0})$
induced by $-E^{*}+\text{Ent}^{*}$ coincides with the \emph{Ding
functional} $\mathcal{D}(u)$ in Kähler geometry, extended to $\mathcal{E}^{1}(X,\omega_{0})$
in \cite{bbgz}. For a general $\beta$ the corresponding twisted
Ding functional $\mathcal{D}_{\beta}$ on $\mathcal{E}^{1}(X,\omega_{0})$
is given by 
\[
\mathcal{D}_{\beta}(u):=-\mathcal{E}(u)+\frac{1}{\beta}\log\int e^{\beta u}dV
\]
An extension of the argument used to prove formula \ref{eq:inf F is inf dual}
(concerning the boundedness statement) now gives
\begin{thm}
\label{thm:M and D}\cite{berm6} The functional $\mathcal{M}_{\beta}$
is bounded from below (coercive) on $\mathcal{E}^{1}(X,\omega_{0})_{0}$
iff $\mathcal{D}_{\beta}$ is bounded from below (coercive) on $\mathcal{E}^{1}(X,\omega_{0})_{0}.$
Moreover, by the regularization result in \cite{bdl1} these properties
are equivalent to the corresponding boundedness/coercivity properties
on the dense subspace $\mathcal{H}(X,\omega_{0})_{0}$ of $\mathcal{E}^{1}(X,\omega_{0})_{0}.$
\end{thm}

For $\beta=-1$ the first statement was first established in \cite{li.h,ru}.
The proof in \cite{li.h} shows that the difference $\mathcal{M}_{\beta}-\mathcal{D}_{\beta}$
is bounded along the Kähler-Ricci flow, thanks to Perelman's estimates,
while the proof in \cite{ru} utilizes the Ricci iteration. In the
case, $\beta=-1$ the coercivity of $\mathcal{M}_{\beta}$ is, in
fact, equivalent to the existence of unique Kähler-Einstein metric,
as first shown in \cite{ti}, using Aubin's method of continuity (discussed
above in connection to Lemma \ref{lem:critical pt}). More recently,
this result has been given a new proof using the notion of geodesics
in $\mathcal{E}^{1}(X)$ and extended in various directions, as discussed
in Section \ref{subsec:The-uniform-YTD}.

\section{\label{sec:The-Yau-Tian-Donaldson-conjectur}The Yau-Tian-Donaldson
conjecture}

\subsection{The Yau-Tian-Donaldson conjecture for polarized manifolds $(X,L)$}

Let $(X,L)$ be a polarized projective algebraic manifold, i.e. $L$
is a holomorphic line bundle over $X$ whose first Chern class $c_{1}(L)$
contains some Kähler form. 
\begin{conjecture}
(Yau-Tian-Donaldson, YTD) There exists a Kähler metric in $c_{1}(L)$
with constant scalar curvature iff $(X,L)$ is K-polystable
\end{conjecture}

We will briefly recall the notion of K-polystability (see the survey
\cite{do2} for more background on the Yau-Tian-Donaldson conjecture
and its relation to Geometric Invariant Theory, GIT). The notion of
K-polystability can be viewed as a ``large $N_{k}-$limit'' of the
classical notion of Chow polystability in GIT with respect to the
action of complex reductive group $GL(N_{k},\C)$ on the Chow variety,
induced from the action of $GL(N_{k},\C)$ on the $N_{k}-$dimensional
complex vector space $H^{0}(X,kL).$ Recall that in GIT the stability
in question is equivalent to the positivity of the GIT-weight of all
one-parameter subgroups (by the Mumford-Hilbert criterion). In the
definition of K-polystability the role of a one-parameter subgroup
$\rho_{k}$ of $GL(N_{k},\C)$ is played by a \emph{test configuration}
$\rho$ for $(X,L).$ In a nutshell, this is a $\C^{*}-$equivariant
embedding 
\[
\rho:(X\times\C^{*},L)\hookrightarrow(\mathcal{X},\mathcal{L})
\]
 of the polarized trivial fibration $(X\times\C^{*},L)$ over $\C^{*}$
into a normal variety $\mathcal{X}$ fibered over $\C$ endowed with
a relatively ample $\Q-$line bundle \emph{$\mathcal{L}.$ }To any
test configuration $\rho$ is attached an invariant, called the \emph{Donaldson-Futaki
invariant} $DF(\rho)\in\R$ and $(X,L)$ is said to be \emph{K-semistable
}if $DF(\rho)\geq0$ for any test configuration,\emph{ K-polystable}
if moreover equality only holds when $\mathcal{X}$ is biholomorphic
to $X\times\C$ and\emph{ K-stable }if the equality only holds when
$\mathcal{X}$ is equivariantly biholomorphic to $X\times\C.$ The
Donaldson-Futaki invariant of $\rho$ may be defined as a limit of
the GIT-weights of a sequence of one-parameter subgroups $\rho_{k}$
of $GL(N_{k},\C)$ induced by $\rho.$ But it may also be expressed
directly as an intersection number \cite{w,od1}:
\[
DF(\rho)=\frac{1}{L^{n}(n+1)}\left(a\mathcal{L}^{n+1}+(n+1)K_{\mathcal{\overline{\mathcal{X}}}/\P^{1}}\cdot\mathcal{L}^{n}\right),\,\,\,a:=-nK_{X}\cdot L^{n-1}/L^{n}
\]
where we have identified a test configuration $(\mathcal{X},\mathcal{L})$
with its $\C^{*}-$equivariant compactification over $\P^{1}$ (obtained
by replacing the base $\C$ of $\mathcal{X}$ with $\P^{1})$ and
the intersection numbers are computed on the compactification $\overline{\mathcal{X}}$
of the total space $\mathcal{X}).$

\subsubsection{\label{subsec:The-uniform-YTD}The uniform YTD and geodesic stability}

The ``only if'' direction of the YTD conjecture was established
in \cite{st} in the case when the group $\text{Aut}(X,L)$ of all
automorphisms of $X$ that lift to $L$ is finite and in \cite{bdl1}
in general. However, for the converse implication there are indications
that the notion of K-polystability needs to be strengthened, in general.
Here we will, for simplicity, focus on the case when $\text{Aut}(X,L)$
is finite. Then K-polystability is equivalent to K-stability and,
moreover, if $c_{1}(L)$ contains a Kähler metric with constant curvature
then it is uniquely determined \cite{do1,ber-bernd}. Following \cite{bhj1,der}
$(X,L)$ is said to be \emph{uniformly K-stable} (in the $L^{1}-$sense)
if there exists $\epsilon>0$ such that
\begin{equation}
DF(\rho)\geq\epsilon\left\Vert \rho\right\Vert _{L^{1}},\label{eq:DG bigger than eps}
\end{equation}
 where the $L^{1}-$norm $\left\Vert \rho\right\Vert _{L^{1}}$ is
defined as the normalized limit of the $l^{1}-$norms of the weights
of the $\C^{*}-$action on the central fiber of $(\mathcal{X},\mathcal{L}).$
The ``only if'' direction of the\emph{ ``uniform YTD conjecture''}
- where K-stability is replaced by uniform K-stability (in the $L^{1}-$sense)
- was established in \cite{bhj1}, by leveraging the connection to
the ``metric space analog'' of the uniform YTD conjecture, to which
we next turn. Denote by $d_{1}$ the metric on $\mathcal{H}(X,\omega_{0})$
induced by the intrinsic $L^{1}-$Finsler metric 
\[
\int_{X}|\dot{u}|^{1}\omega_{u_{0}}^{n},\,\,\,\,\dot{u}:=\frac{du}{dt}|_{t=0},\,\,\,u_{0}\in\mathcal{H}.
\]
 As shown in \cite{da} the metric space completion $\left(\overline{\mathcal{H}(X,\omega_{0})_{0}},d_{1}\right)$
may be identified with the space $\mathcal{E}^{1}(X,\omega_{0})_{0}$
(discussed in Section \ref{subsec:Pluripotential-theory-recap}) and
$d_{1}(u,0)$ is comparable to $-\mathcal{E}(u),$ which equivalently
means that there exists a constant $c$ such that

\begin{equation}
-c+c^{-1}d_{1}(u,0)\leq E\left(MA(u)\right)\leq cd_{1}(u,0)+c.\label{eq:d one compared to E}
\end{equation}
The relevant constant speed geodesics $u_{t}$ in the metric space
$\left(\mathcal{E}^{1}(X,\omega_{0})_{0},d_{1}\right)$ have the property
that 
\begin{equation}
U(x,\tau):=u_{-\log|\tau|}(x)\in PSH(X\times D^{*},\omega_{0}),\label{eq:def of U}
\end{equation}
 where we are using the same notation $\omega_{0}$ for the pull-back
of $\omega_{0}$ to the product $X\times D^{*}$ of $X$ with the
punctured unit-disc $D^{*}$ in $\C.$ In fact, $u_{t}$ may be characterized
by a maximality property of the corresponding $\omega_{0}-$psh function
$U$ \cite{bbj}. Any test configuration $\rho$ induces a geodesic
ray $u_{t}$ in $\mathcal{E}^{1}(X,\omega_{0})_{0},$ emanating from
$0\in\mathcal{H}(X,\omega_{0})$ (such that $U$ extends, after removing
divisorial singularities, to a bounded function on $\mathcal{X}$)
\cite{p-s,da}. Moreover, 
\[
\left\Vert \rho\right\Vert _{L^{1}}=\frac{d}{dt}d_{1}(u_{t},0)=t^{-1}d(u_{t},0)
\]
 for any $t>0.$ As conjectured in \cite{ch0}, and confirmed in \cite{ber-bernd},
the Mabuchi functional $\mathcal{M}$ (Section \ref{subsec:The-Mabuchi-and})
is convex along geodesic $u_{t}$ such that $\omega_{U}\in L_{loc}^{\infty}.$
More generally, the extension of $\mathcal{M}$ to $\mathcal{E}^{1}(X,\omega_{0})$
is also convex along geodesics $u_{t}$ \cite{bdl1}. In particular,
its (asymptotic) slope 
\[
\mathcal{\dot{M}}(u_{t}):=\lim_{t\rightarrow\infty}t^{-1}\mathcal{M}(t)\in]-\infty,\infty]
\]
 is well-defined. In the case when $u_{t}$ is the geodesic ray attached
to a test configuration $\rho$ the slope $\mathcal{\dot{M}}(u_{t})$
is closely related to $DF(\rho)$ (the two invariants coincide after
a base change \cite{sj-d,li2}). 
\begin{thm}
\label{thm:metric ytd} \cite{d-r,b-d-l2,c-cII} Let $(X,L)$ be a
polarized manifold. The following is equivalent.
\begin{enumerate}
\item $(X,L)$ admits a unique Kähler metric with constant scalar curvature 
\item $(X,L)$ is geodesically stable, i.e. $\mathcal{\dot{M}}(u_{t})>0$
for any non-trivial geodesic ray $u_{t}$ in $\mathcal{E}^{1}(X,\omega_{0})_{0}$
\item $\mathcal{M}$ is coercive on $\mathcal{E}^{1}(X,\omega_{0})_{0}$
(or, equivalently, on $\mathcal{H}(X,\omega_{0})_{0}\subset\mathcal{E}^{1}(X,\omega_{0})_{0})$
\end{enumerate}
\end{thm}

The equivalence $"2\iff3$'' is implicit in \cite{d-r} (see \cite[Thm 2.16]{bbj}
for a generalization). It can be seen as an analog of the classical
fact that a convex function on Euclidean $\R^{n}$ is comparable to
the distance to the origin iff all its slopes are positive. In the
proof of $"2\iff3$'' a substitute for the compactness of the unit-sphere
in $\R^{n}$ (parametrizing all unit speed geodesics) is provided
by the energy-entropy compactness in Theorem \ref{thm:(energy/entrop-compactness)}.
The implication $"1\implies3"$ follows directly from the convexity
of $\mathcal{M}$ combined with the weak-strong uniqueness result
in \cite{b-d-l2}, showing, in particular, that if $(X,L)$ admits
a unique Kähler metric with constant scalar curvature $\omega,$ then
any minimizer of $\mathcal{M}$ in $\mathcal{E}^{1}$ coincides with
the Kähler potential of $\omega.$ The final implication $"3\implies1\text{\textquotedblright}$
was recently settled in \cite{c-cII}, using a new a priori estimate
for a generalization of Aubin's continuity method for constant scalar
curvature metrics (bounding the $C^{0}-$norm of the solutions by
the entropy of the corresponding Monge-Ampère measures, which, in
turn is uniformly bounded under the coercivity assumption).

\subsection{\label{subsec:The-variational-approach}The variational approach
to the uniform YTD conjecture in the ``Fano case''}

The ``Fano case'' of the YTD conjecture, i.e. the case when $X$
is Fano and $L=-K_{X},$ was settled in \cite{c-d-s}, by establishing
Tian's partial $C^{0}-$estimate \cite{ti0b} along a singular version
of Aubin's continuity method. Here we will focus on the variational
proof of the uniform YTD conjecture on Fano manifolds in \cite{bbj},
which, in particular, exploits the notion of Ding stability originating
in \cite{berman6ii} (as further developed in \cite{bhj1,bbj}; see
the survey \cite{bo} for more background). 
\begin{thm}
\cite{bbj} \label{thm:variational ytd} Let $X$ be a Fano manifold.
The following is equivalent:
\begin{enumerate}
\item $X$ admits a unique Kähler-Einstein-metric
\item $X$ is uniformly Ding stable
\item $X$ is uniformly K-stable 
\end{enumerate}
\end{thm}

The implication $"1\implies2"$ follows from the convexity of the
Ding functional along geodesics, as in \cite{berman6ii} - here we
shall focus on the converse implication. By Theorem \ref{thm:metric ytd}
it is enough to show that if $X$ is uniformly Ding stable, then $X$
it geodesically stable. This is achieved in \cite{bbj}, using a valuative
(non-Archimedean) language. For simplicity, it may be helpful to briefly
first describe the argument with the non-Archimedean language stripped
away. The starting point is the observation that the function $U$
on $X\times D^{*}$ corresponding to a geodesic $u_{t}$ in $\mathcal{E}^{1}(X,\omega_{0})_{0}$
(formula \ref{eq:def of U}) extends to a sup-normalized $\omega_{0}-$psh
function $U$ on $X\times D,$ which, however, is highly singular
on $X\times\{0\},$ unless $u_{t}$ is trivial. But employing Demailly's
approximation procedure \cite{dem} (involving the multiplier ideal
sheaves $\mathfrak{J}(kU),$whose definition is recalled in the following
section) the function $U$ may be expressed as a decreasing limit
of $S^{1}-$invariant $\omega_{0}$-psh functions $U_{k}$ with analytic
(algebraic) singularities, which define $\C^{*}-$invariant ideals
$\mathfrak{J}_{k}$ supported in $X\times\{0\}.$ Accordingly, by
standard resolution of singularities there exists a $\C^{*}-$equivariant
holomorphic surjection $\pi_{k}$ from a non-singular variety $\mathcal{X}_{k}$
to $X\times\C$ such that $E_{k}:=\pi_{k}^{*}\mathfrak{J}_{k}$ is
a principal ideal, i.e. defines a divisor on $\mathcal{X}_{k}.$ This
procedure yields a sequence of test-configurations $\rho_{k}=(\mathcal{X}_{k},\mathcal{L}_{k})$
where $\mathcal{L}_{k}$ is the pull-back to $\mathcal{X}_{k}$ of
$L\rightarrow X$ with an appropriate multiple of $\mathcal{O}(E_{k})$
subtracted. To show that ``$3\implies1"$ it would, essentially,
be enough show that the slope $\mathcal{M}(u_{t})$ dominates the
Donaldson-Futaki invariants $DF(\mathcal{\rho}_{k}).$ However, this
leads to technical problems that are bypassed by exploiting that $\mathcal{M}\geq\mathcal{D},$
where $\mathcal{D}$ is the Ding functional on \emph{$\mathcal{H}_{0}$}
(discussed in Section \ref{subsec:The-Mabuchi-and}) which behaves
better under the approximation procedure above, giving 
\begin{equation}
\mathcal{\dot{D}}(u_{t})\geq\liminf_{k\rightarrow\infty}\mathcal{D}(\rho_{k}),\label{eq:D slope bigger than D inv}
\end{equation}
 where $\mathcal{D}(\mathcal{\rho}_{k})$ is the ``Ding invariant''
originating in \cite{berman6ii} (that we shall come back to in Section
\ref{subsec:The-thermodynamical-formalism NA}). Assuming that $\mathcal{X}$
is uniformly Ding stable this shows that $"2\implies1"$ (after a
twist of the argument which amounts to replacing $\mathcal{D}$ with
$\mathcal{D}_{\beta}$ for $\beta=-(1+\epsilon)$ ).

Finally, the equivalence $"2\iff3"$ is shown in the first preprint
version of \cite{bbj}, using techniques from the Minimal Model Program,
inspired by \cite{l-x} (the proof can - loosely speaking - be interpreted
as a non-Archimedean analog of the Kähler-Ricci flow argument in \cite{li.h}
mentioned in connection to Theorem \ref{thm:M and D}). The equivalence
$"2\iff3"$ in the general setup of log Fano varieties is established
in \cite{fu2}.

\subsubsection{Twisted Kähler-Einstein metrics}

The results in \cite{bbj} apply more generally to Kähler-Einstein
metrics twisted by a positive klt current $\theta,$ showing that
such a metric exists iff $\delta_{\theta}(X)>1,$ where $\delta_{\theta}(X)$
is a twisted generalization of the invariant $\delta(X)$ appearing
in formula \ref{eq:def of delta as limit of delta k}. This part of
the proof does not need any results from the Minimal Model Program
(as discussed in the following section). As a corollary it is also
shown that
\begin{equation}
\min\{1,\delta(X)\}=\min\{1,\Gamma(X)\}=R(X)\label{eq:min delta X}
\end{equation}
 where $\Gamma(X)$ is the invariant appearing in Conjecture \ref{conj:lim lct is Gamma}
and $R(X)$ denotes the greatest lower bound on the Ricci curvature
(independently shown in \cite{c-r-z}).

\subsection{\label{subsec:Non-Archimedean-pluripotential-t}Non-Archimedean pluripotential
theory and the variational formula for $\delta(X)$}

The only properties of the geodesic $u_{t}$ that actually entered
into the proof outlined above concerned the multiplier ideal sheaves
$\mathfrak{J}(kU)$ of the $\omega_{0}-$psh function $U$ on $X\times D,$
whose stalks consist of all germs of holomorphic functions $f$ such
that $|f|^{2}e^{-2kU}$ is locally integrable. In turn, the multiplier
ideal sheaves $\mathfrak{J}(kU)$ only depend on the Lelong numbers
of $U$ on all modifications (blow-ups) of $X\times\C$ (see \cite[Thm A]{bfj}
and Thm B.5 in Appendix B of \cite{bbj}). The Lelong numbers in question
can be packaged into a function $U(v)$ on the space $[X\times\C]_{\text{div}}$
of all divisorial valuations $v$ on $X\times\C,$ as follows. First
recall that, by definition, a divisorial valuation $v$ on variety
$Y$ is encoded by a positive number $c$ and a prime divisor $E_{v}$
over $Y,$ i.e a prime divisor on some blow-up of $Y$ (which may
be assumed to be a non-singular hypersurface). Such a valuation $v$
acts on rational (meromorphic) function $f\in\C(Y)$ by $v(f):=c\text{ord}_{E_{v}}(f)\in\R,$
where $\text{ord}_{E_{v}}(f)$ denotes the order of vanishing at a
generic point of $E_{v}$ of the pull-back of $f.$ Now, if $U$ is,
locally, of the form $U=\log|f|+O(1)$ for a holomorphic function,
one defines
\[
U(v):=-v(f):=-c\text{ord}_{E_{v}}(f).
\]
 In the general definition of $U(v)$ one replaces $\text{ord}_{E_{v}}(f)$
with the Lelong number of $U$ at a generic point $p$ of $E_{v}$
(i.e. the sup of all $\lambda\in[0,\infty[$ such that $f\leq\lambda\log|z|+O(1)$
wrt local holomorphic coordinates $z$ centered at $p).$ In this
context Demailly's approximation procedure yields
\begin{equation}
U_{k}(v):=k^{-1}\max_{i}\left(-\text{ord}_{E_{v}}(f_{i}^{(k)})\right)\rightarrow U(v),\label{eq:def of U k}
\end{equation}
 where $f_{i}^{(k)}$ denote local generators of the multiplier ideal
sheaf $\mathfrak{J}(kU).$ In fact, after passing to a subsequence
(replacing $k$ with $2^{k})$ the sequence $U_{k}$ is decreasing
in $k$ (by the subaditivity of multiplier ideals). 

\subsubsection{Pluripotential theory on the Berkovich space $X_{NA}$}

In the present setup the valuative procedure above is initially applied
to $Y=X\times\C.$ However, exploiting that we are only interested
in the value $U(w)$ at a divisorial valuation $w$ on $X\times\C$
which is $\C^{*}-$invariant, we can identify $[U](w)$ with the function
on $u(v)$ on $X_{\text{div}},$defined by 
\[
u(v):=U(w),\,\,\,\,v\in X_{\text{div}},\,\,w\in(X\times\C)_{\text{div}}
\]
 where $w$ is the Gauss extension of $v,$ defining a $\C^{*}-$equivariant
valuation over $X\times\C$ normalized by $w(\tau)=1$ (where $\tau$
denotes the coordinate on the factor $\C)$ \cite[Section 4.1]{bhj1}.
Next, by identifying a valuation $v$ on $X$ with the corresponding
non-Archimedean absolute value on $\C(X),$ i.e. with $|\cdot|_{v}:=e^{-v(\cdot)},$
the space $X_{\text{div}}$ injects as a dense subspace of the Berkovich
analytification $X_{NA}$ of the projective variety $X$ over the
field $\C,$ induced by the trivially valued absolute value on the
ground field $\C$ (locally consisting of all multiplicative semi-norms
extending the trivially valued absolute value, $|\cdot|_{v}\equiv1,$
on the field $\C$). The notation $X_{NA}$ (with $NA$ a shorthand
for Non-Archimedean) is used here to distinguish $X_{NA}$ from $X$
which is the Berkovich analytification in the ``Archimedean case'',
i.e. the case of the standard absolute value $|\cdot|$ on the ground
field $\C.$ 

The topological space $X_{NA}$ has the virtue of being both compact
and connected. Moreover, the function $u(v)$ on $X_{\text{div}}$
extends to a plurisubharmonic (psh) function on $X_{NA}$ in the sense
of \cite{b-j0}, denoted by $u_{NA}.$ Indeed, in analogy to the Archimedean
case one can first define $\mathcal{H}(X_{NA})_{0}$ to be the space
of all functions $u_{NA}$ on $X_{NA}$ induced by test-configurations
$\rho,$ as above, and then define $PSH(X_{NA})$ as the space of
all functions that can be written as decreasing nets of functions
in $\mathcal{H}(X_{NA})_{0}$ plus constants (functions in $PSH(X_{NA})$
are called $L-$psh in \cite{b-j0} to emphasize their global dependence
on $L$). There is a Monge-Ampère operator $MA$ on $\mathcal{H}(X_{NA})$
taking values in the space of probability measures on $X_{NA}$ \cite{bhj1,b-j0}
(which, in a very general setup can be defined in terms of the non-Archimedean
generalization of exterior products of curvature forms introduced
in \cite{c-d}). Concretely, $MA(u_{NA})$ is a discrete probability
measure supported on the valuations $v_{i}\in X_{\text{div}}$ induced
by irreducible components of the central fiber of the test configuration
corresponding to $u_{NA}$ \cite[Section 6.7]{bhj1}. Anyhow, in the
present setup one may directly define $MA$ on $\mathcal{H}(X_{NA})$
as the differential of the functional 
\[
\mathcal{E}_{NA}(u_{NA}):=\frac{\mathcal{L}^{n+1}}{(n+1)L^{n}},
\]
 whose definition mimics formula \ref{eq:formula for e u} (with $\omega_{0}=0);$
this analogy becomes more clear when both $\mathcal{E}$ and $\mathcal{E}_{NA}$
are expressed in terms of Deligne pairings \cite{b-e}. As in the
usual Archimedean setup (Section \ref{subsec:Pluripotential-theory-recap})
the function $\mathcal{E}_{NA}$ on $\mathcal{H}(X_{NA})$ has a unique
smallest usc extension to $PSH(X_{NA});$ the subspace $\left\{ \mathcal{E}_{NA}>-\infty\right\} $
of $PSH(X_{NA})$ is denoted by $\mathcal{E}^{1}(X_{NA})$ and $MA$
extends to $\mathcal{E}^{1}(X_{NA}),$ as the differential of the
functional $\mathcal{E}_{NA}^{1}.$
\begin{rem}
\label{rem:maximal rays}The map $u_{t}\mapsto u_{NA}$ from geodesic
rays in $\mathcal{E}^{1}(X,\omega_{0})_{0}$ to the space $\mathcal{E}^{1}(X_{NA})_{0},$
described above, has the property that $\dot{\mathcal{E}}(u_{t})\leq\mathcal{E}(u_{NA})$
and is, in general, not injective. The geodesic rays satisfying $\dot{\mathcal{E}}(u_{t})=\mathcal{E}(u_{NA})$
are precise those called \emph{maximal} in \cite[Section 6.4]{bbj}
and they are in one-to-one correspondence with $\mathcal{E}^{1}(X_{NA}).$
\end{rem}

\subsubsection{\label{subsec:The-thermodynamical-formalism NA}The thermodynamical
formalism}

The non-Archimedean formalism naturally ties in with the thermodynamical
formalism (discussed in Section \ref{sec:The-thermodynamical-formalism}).
For example, as shown \cite{bhj1,b-j,b-j0}, up to a base change of
$\rho,$\footnote{the base change is needed as the rhs in formula \ref{eq:DF as NA MAb}
is one-homogenuous under the natural action of $\R_{>0}$ on $X_{NA},$
corresponding to a base change of $\rho.$} 
\begin{equation}
DF(\rho)=\mathcal{M}_{NA}(U_{NA}):=F_{NA}(MA(U_{NA})),\label{eq:DF as NA MAb}
\end{equation}
 where $F_{NA}$ is the non-Archimedean analog on $\mathcal{P}(X_{NA})$
of the free energy functional $F$ on $\mathcal{P}(X)$ defined by
\[
F_{NA}(\mu)=-E_{NA}(\mu)+\text{Ent}_{NA}(\mu),
\]
 where the non-Archimedean energy $E_{NA}(\mu)$ may be defined as
a Legendre-Fenchel transform of the functional $\mathcal{E}_{NA}$
and the non-Archimedean entropy $\text{Ent}_{NA}(\mu)$ is defined
by
\[
\text{Ent}_{NA}(\mu):=\int_{X_{NA}}A(v)\mu,\,\,\,\,A(v):=c\left(1+\text{ord}_{E_{v}}(K_{Y_{v}/X})\right)\,v\in X_{\text{div}}
\]
 where $A(v)$ is the\emph{ log discrepancy}, defined as the greatest
lsc extension to $X_{NA}$ of the function on $X_{\text{div}}$ defined
above. Thus, in contrast to the usual entropy functional on $\mathcal{P}(X),$
the non-Archimedean entropy is a linear functional. Likewise, the
``Ding invariant'' appearing in formula \ref{eq:D slope bigger than D inv}
may be expressed as follows in terms of the Legendre-Fenchel transform
\[
\mathcal{D}(\rho)=\mathcal{D}_{NA}(u_{NA}):=-E_{NA}^{*}(u_{NA})+\text{Ent}_{NA}^{*}(u_{NA})
\]
 in analogy with the usual Archimedean setup in Section \ref{subsec:The-Mabuchi-and}.
The inequality \ref{eq:D slope bigger than D inv} is then obstained
by showing that the slope $\mathcal{\dot{D}}(u_{t})$ is bounded from
below by $\mathcal{D}(u_{NA}),$ which, in turn, equals the limit
of $\mathcal{D}(\rho_{k})$ (where $\rho_{k}$ is the test configuration
corresponding to $U_{k}$ defined by formula \ref{eq:def of U k}). 

As shown in \cite{b-j} (and \cite{bbj} in the general twisted setting)
the thermodynamical formalism can be leveraged to prove the following
theorem ($"1\iff3"$ is shown in \cite{fu2} using the Minimal Model
Program):
\begin{thm}
\label{thm:-delta iff unif K iff Ding}\cite{b-j} Let $X$ be a Fano
manifold. The following is equivalent: 
\begin{enumerate}
\item $\delta(X)>1$
\item $X$ is uniformly K-stable on $\mathcal{E}^{1}(X_{NA})$ (i.e. the
inequality \ref{eq:DG bigger than eps} extends from $\mathcal{H}(X_{NA})$
to $\mathcal{E}^{1}(X_{NA})$)
\item $X$ is uniformly Ding stable
\end{enumerate}
\end{thm}

The starting point of the proof of $"1\iff2"$ is the following variational
formula for $\delta(X)$ established in \cite{bl-j,b-j}, realizing
$\delta(X)$ as a ``stability threshold'' (where $\delta_{v}$ denotes
the Dirac measure at a point $v$ in $X_{NA}):$ 
\begin{equation}
\delta(X)=\inf_{v\in X_{\text{div}}}\frac{\text{Ent}_{NA}(\delta_{v})}{E_{NA}(\delta_{v})}=\inf_{v\in X_{NA}}\frac{\text{Ent}_{NA}(\delta_{v})}{E_{NA}(\delta_{v})}=\inf_{\mu\in\mathcal{P}(X_{NA})}\frac{\text{Ent}_{NA}(\mu)}{E_{NA}(\mu)}\label{eq:delta equal to Ent over E}
\end{equation}
 using, in the second equality, that $X_{\text{div}}$ is dense in
$X_{NA}$ (together with a semi-continuity argument) and in the last
equality (shown in \cite{b-j}) that $\text{Ent}_{NA}(\mu)$ and $E_{NA}(\mu)$
are linear and convex, respectively, on $\mathcal{P}(X_{NA}).$ The
function $v\mapsto E_{NA}(\delta_{v})$ is usually denoted by $S(v)$
and can can be shown to coincide with the ``expected order of vanishing
along $v"$ \cite{bl-j}. In terms of  the non-Archimedean version
of the free energy functional at inverse temperature $\beta,$ denoted
by $F_{NA,\beta}(\mu),$ formula \ref{eq:delta equal to Ent over E}
yields
\[
\delta(X)\geq1+\epsilon\iff\inf_{\mu\in\mathcal{P}(X_{NA})}F_{NA,-1-\epsilon}(\mu)\geq0\iff\inf_{\mu\in\mathcal{P}(X_{NA})}\frac{F_{NA}(\mu)}{E_{NA}(\mu)}\geq\epsilon.
\]
 Finally, expressing $\mu=MA(U_{NA})$ for $U_{NA}\in\mathcal{E}^{1}(X_{NA}),$
using the non-Archimedean version of the ``Calabi-Yau correspondence''
\ref{eq:u maps to MA}, and invoking the non-Archimedean version of
the inequalities \ref{eq:d one compared to E} (established in \cite{bhj1})
proves the equivalence $"1\iff2".$ Next, using the Legendre-Fenchel
transform, just as in the proof of Theorem \ref{thm:M and D}, one
sees that uniform K-stability on $\mathcal{E}^{1}(X_{NA})$ is equivalent
to uniform Ding stability on $\mathcal{E}^{1}(X_{NA}).$ Finally,
$"2\iff3"$ follows from the fact that $\mathcal{D}_{NA}$ is continuous
under approximation of $U_{NA}\in\mathcal{E}^{1}(X_{NA})$ by a decreasing
sequence in $\mathcal{H}(X_{NA})$ (e.g. using multiplier ideal sheaves
as in formula \ref{eq:def of U k}).

In order to deduce the equivalence $"2\iff3"$ in Theorem \ref{thm:variational ytd}
from the previous theorem it would be enough to prove the following
non-Archimedean analog of the regularization property shown in \cite[Section 3]{bdl1}.
\begin{conjecture}
\label{conj:reg}\cite{b-j} Given any $u\in\mathcal{E}^{1}(X_{NA})$
there exists a sequence of $u_{j}\in\mathcal{H}(X_{NA})$ converging
weakly towards $u$ such that $E_{NA}\left(MA(u_{j})\right)$ and
$\text{Ent}_{NA}\left(MA(u_{j})\right)$ converge towards $E_{NA}\left(MA(u)\right)$
and $\text{Ent}_{NA}\left(MA(u)\right),$ respectively.
\end{conjecture}

\subsection{Recent developments }

Recently there has been an explosion of exciting further developments.
In \cite{ltw,li1} Theorem \ref{thm:variational ytd} and its variational
proof was extended to general singular (log) Fano varieties using,
in particular, the singular version of Theorem \ref{thm:reg and exist of min}
established in \cite{bbegz}. Moreover, very recently it was shown
in \cite{l-x-z}, using techniques from the Minimal Model Program,
that the infimum over $X_{\text{div}}$ in formula \ref{eq:delta equal to Ent over E}
is (when $\delta(X)\leq1)$ attained at some $v\in X_{\text{div}}.$
Moreover, any such minimizing divisorial valuation $v$ has the property
that associated graded ring is finitely generated and defines a special
test configuration $\rho$ for $(X,-K_{X}).$ In particular the central
fiber of $\rho$ is irreducible (the relation between test configurations,
filtrations and finitely graded rings originates in \cite{wn,sz}).
In non-Archimedean terms the result in \cite{l-x-z} can be formulated
as a regularity result for the minimizer in question, saying that
$\delta_{v}=MA(U_{NA})$ for some $U_{NA}\in\mathcal{H}(X_{NA})$
(in analogy to the regularity result in Theorem \ref{thm:reg and exist of min};
cf. the appendix in \cite{fu2}). As a corollary it is shown in \cite{l-x-z}
that uniform K-stability is equivalent to K-stability. In fact, these
results are shown to hold in the general setup of (log) Fano varieties.
When combined with the aforementioned results in \cite{ltw,li1} this
settles the YTD conjecture in the general setting of (log) Fano varieties
(the ``only if'' implication was previously shown in \cite{berman6ii}).
In another direction, a new variational proof of the uniform YTD conjecture
in the non-singular Fano case is given in \cite{zh}, using the quantized
Ding-functional (leveraging the result in \cite{rtz} saying that
the algebro-geometric invariant $\delta_{k}(X)$ in formula\ref{eq:def of delta k}
concides with coercivity threshold of the quantized Ding-functional).
More generally, the results in \cite{zh} imply that the first equality
in formula \ref{eq:min delta X} holds without taking the minimum
with $1$ (by combining \cite{zh} with Theorem \ref{thm:M and D})

The variational/non-Archimedean approach is extended to polarized
manifolds $(X,L)$ in \cite{li2} to show that, if $X$ is uniformly
K-stable on $\mathcal{E}^{1}(X_{NA})$ (as in Theorem \ref{thm:-delta iff unif K iff Ding}),
then $X$ is geodesically stable and thus by Theorem \ref{thm:metric ytd}
(i.e. by \cite{c-cII}) $(X,L)$ admits a Kähler metric with constant
scalar curvature. The converse statement is, however, still open.
The complete solution of the uniform YTD conjecture for $(X,L)$ is
thus reduced to Conjecture \ref{conj:reg}. An important ingredient
in \cite{li2} is the notion of maximal geodesic rays $u_{t}$ introduced
in \cite{bbj} (see Remark \ref{rem:maximal rays}). The theory of
maximal geodesic rays is further developed in in \cite{d-x} and related
to singularity types of quasi-psh functions and the Legendre transform
construction of geodesic rays introduced in \cite{r-w}. In \cite{x}
analytic variants of stability thresholds are introduced, expressed
in terms of singularity types of quasi-psh functions.

\section{\label{sec:A-Non-Archimedean-approach}A Non-Archimedean approach
to Gibbs stability}

This final section is a report on joint work in progress with Sébastien
Boucksom and Mattias Jonsson to prove the converse of Theorem \ref{thm:(Fujita-Odaka)-Uniform-Gibbs}
or, more generally, to prove that
\begin{equation}
\lim_{N\rightarrow\infty}\text{lct }(D_{N})=\delta(X)\label{eq:lim D N is delta}
\end{equation}
 (which, when combined with results in \cite{zh}, would also settle
Conjecture \ref{conj:lim lct is Gamma}). The strategy is to adapt
the variational approach to the convergence in Conjecture \ref{conj:Fano with triv autom intr},
discussed in Section \ref{subsec:The-case beta negative proof}, to
the non-Archimedean setup. The starting point is the standard valuative
expression for the log canonical threshold of a divisor that yields
(using the notation in Section \ref{subsec:Non-Archimedean-pluripotential-t})

\begin{equation}
\text{lct }(D_{N})=\inf_{v^{(N)}\in[X^{N}]_{\text{div}}}\frac{A(v^{(N)})}{k^{-1}\left(v^{(N)}(\det S^{(k)})\right)}:=\frac{N^{-1}A(v^{(N)})}{E_{NA}^{(N)}(v^{(N)}},\label{eq:lct as inf}
\end{equation}
 where we have introduced the \emph{non-Archimedean energy per particle}
as the following function on $[X^{N}]_{\text{div}}:$
\[
E_{NA}^{(N)}(v^{(N)}):=N^{-1}k^{-1}\left(v^{(N)}(\det S^{(k)})\right)=:-N^{-1}k^{-1}\log\left|\det S^{(k)}\right|_{v^{(N)}}
\]
 (which is proportional to the negative of the psh function on $[X^{N}]_{NA}$
induced by the quasi-psh function $\log\left\Vert \det S^{(k)}\right\Vert ^{2}$
on $X^{N}).$ In this notation formula \ref{eq:lct as inf} can be
viewed as a non-Archimedean analog of Gibbs variational principle
\ref{eq:Gibbs var pr} (since $\text{lct }(D_{N})-1$ is equal to
the one-homogeneous non-Archimedean ``N-particle free energy'' $-E_{NA}^{(N)}+N^{-1}A,$
normalized by $E_{NA}^{(N)}).$ There are standard inclusions $i_{N}$
and surjections $\pi_{N}$, 
\[
i_{N}:\,\left(X_{NA}\right)^{N}\hookrightarrow[X^{N}]_{NA},\,\,\,\pi_{N}:\,[X^{N}]_{NA}\twoheadrightarrow\left(X_{NA}\right)^{N}.
\]
(the map $i_{N}$ is, however, not surjective). The non-Archimedean
version of the empirical measure $\delta_{N}$ mapping $\left(X_{NA}\right)^{N}$
to $\mathcal{P}(X_{NA})$ (obtained by replacing $X$ with $X_{NA}$
in formula \ref{eq:def of empirical measure}) thus induces a map
\[
\pi_{N}^{*}\delta_{N}:\,\,[X^{N}]_{\text{div}}\rightarrow\mathcal{P}(X{}_{NA}),\,\,\,v^{(N)}\mapsto N^{-1}\sum_{i=1}^{N}\delta_{\left(\pi_{N}(v^{(N)})\right)_{i}}.
\]
 It follows from the results in \cite{zh} (which are non-Archimedean
versions of results in \cite{b-b}) that the restriction of $E_{NA}^{(N)}$
to $(X_{NA})^{N}$ Gamma-converges towards $E_{NA}(\mu)$ (in analogy
with the convergence \ref{eq:Gamma conv}). In particular,
\begin{equation}
\lim_{N\rightarrow\infty}\delta_{N}(v_{1},...,v_{N})=\mu\in\mathcal{P}([X]_{NA})\implies\liminf_{N\rightarrow\infty}E_{NA}^{(N)}\left(i_{N}(v_{1},...,v_{N})\right)\geq E_{NA}(\mu)\label{eq:lower bound on N particle NA energ}
\end{equation}
Moreover, $N^{-1}A\left(i_{N}(v_{1},...,v_{N})\right)=\int_{X_{NA}}A(v)\delta_{N}(v_{1},...,v_{N}),$
as follows readily from the definitions. Hence, restricting the inf
in formula \ref{eq:lct as inf} to $v_{N}$ of the form $v_{N}=i_{N}(v,...,v)$
for $c\in X_{div}$ reveals that the limsup of $\text{lct }(D_{N})$
is bounded from above by $A(v)/E(\delta_{v})$, proving the upper
bound in formula \ref{eq:lim D N is delta}. This proof essentially
amounts to a reformulation of the proof of Theorem \ref{thm:(Fujita-Odaka)-Uniform-Gibbs}
in \cite{f-o} into a non-Archimedean language. But the main point
of the non-Archimedean formulation is that it opens the door for a
non-Archimedean approach to the missing lower bound. Indeed, it can
be shown that 
\[
\lim_{N_{j}\rightarrow\infty}(\pi_{N_{j}}^{*}\delta_{N_{j}})(v^{(N_{j})})=\mu\in\mathcal{P}\left([X]_{NA}\right)\implies\liminf_{N_{j}\rightarrow\infty}N_{j}^{-1}A(v_{N_{j}})\geq\text{Ent}_{NA}(\mu)
\]
 Hence, all that remains is to establish the following hypothesis
for any valuation $v_{*}^{(N)}$ realizing the infimum in formula
\ref{eq:lct as inf} (which is a non-Archimedean analog of the hypothesis
\ref{eq:hyp for mean en}):
\begin{equation}
\text{Hypothesis: \ensuremath{\,\,\,\lim_{N_{j}\rightarrow\infty}(\pi_{N_{j}}^{*}\delta_{N_{j}})(v_{*}^{(N_{j})})=\mu_{*}\in\mathcal{P}(X_{NA})\implies}\ensuremath{\limsup_{N_{j}\rightarrow\infty}}}E_{NA}^{(N_{j})}(v_{*}^{(N_{j})})\leq E_{NA}(\mu_{*}),\label{eq:hypo}
\end{equation}
(by \ref{eq:lower bound on N particle NA energ} the opposite inequality
holds). Indeed, if the hypothesis holds then we get
\begin{equation}
\inf_{\mu\in\mathcal{P}([X]_{NA})}\frac{\text{Ent}_{NA}(\mu)}{E_{NA}(\mu)}\leq\liminf_{N\rightarrow\infty}\text{lct }(D_{N})\leq\limsup_{N\rightarrow\infty}\text{lct }(D_{N})\leq\inf_{v\in[X]_{div}}\frac{\text{Ent}_{NA}(\delta_{v})}{E_{NA}(\delta_{v})},\label{eq:inequalities in the end}
\end{equation}
which, when combined with the identity \ref{eq:delta equal to Ent over E},
yields the desired formula \ref{eq:lim D N is delta}.

It remains to verify the inequality in the hypothesis above. It would
be enough to establish the following \emph{``restriction hypothesis'':}
the minimizer $v_{*}^{(N)}$ can, asymptotically, be taken to be of
the form $i_{N}(v_{*},v_{*},...,v_{*})$ for a fixed divisorial valuation
$v_{*}$ on $X,$ i.e. 
\[
\exists v_{*}\in X_{\text{div }}\text{such that}\,\,\liminf_{N\rightarrow\infty}\text{lct }(D_{N})=\liminf_{N\rightarrow\infty}\frac{N^{-1}A\left(i_{N}(v_{*},v_{*},...,v_{*})\right)}{E_{NA}^{(N)}\left(i_{N}(v_{*},v_{*},...,v_{*})\right)}.
\]
Indeed, it follows from the convergence of Fekete points on $X_{NA}$
in \cite{b-e} that 
\begin{equation}
\lim_{N\rightarrow\infty}E_{NA}^{(N)}\left(i_{N}(v,v,...,v)\right)=E(\delta_{v})\label{eq:limit of NA E}
\end{equation}
for \emph{any} divisorial valuation $v$ on $X$ (or more generally:
for any non-pluripolar point $v$ in $X_{NA}).$ In particular, it
then follows that any $v_{*}$ satisfying the ``restriction hypothesis''
above computes $\delta(X).$ For instance, it can be verified that
the ``restriction hypothesis'' does hold for log Fano curves $(X,\Delta).$
Anyhow, for any given divisorial valuation $v$ on $X$ formula \ref{eq:limit of NA E}
yields a ``microscopic'' formula for the non-Archimedean free-energy
$F_{NA}(\delta_{v})$ (coinciding with the invariant $\beta(v)$ introduced
in \cite{fu2}) of independent interest:
\[
F_{NA}(\delta_{v}):=-E(\delta_{v})+A(\delta_{v})=\lim_{N\rightarrow\infty}\left(-E_{NA}^{(N)}\left(i_{N}(v,v,...,v)\right)+N^{-1}A\left(i_{N}(v,v,...,v)\right)\right).
\]
 In particular, if $\rho$ is a given test configuration, whose central
fiber $\mathcal{X}_{0}$ is irreducible, this gives a new formula
for the Donaldson-Futaki invariant $DF(\rho),$ using that $DF(\rho)=F_{NA}(\delta_{v}),$
where $v$ is the divisorial valuation on $X$ corresponding to $\mathcal{X}_{0}.$
Comparing with the formula for $DF(\rho)$ in terms of Chow weights
thus suggests that the divisorial valuation $i_{N}(v,v,...,v)$ on
$X^{N},$ attached to $v,$ plays the role of the one-parameter subgroup
of $GL(N,\C)$ attached to $\rho.$ Accordingly, the ``restriction
hypothesis'' is an analog of the Hilbert-Mumford criterion for stability
in Geometric Invariant Theory. 

Finally, coming back to the statistical mechanical point of view discussed
in Section \ref{subsec:The-statistical-mechanical} it may be illuminating
to point out that the ``restriction hypothesis'' essentially amounts
to a concentration phenomenon which may be pictured as follows. Let
us decrease the inverse temperature $\beta$ from a given positive
value towards the critical negative inverse temperature $\beta_{N}$
where $\mathcal{Z}_{N}(\beta)=\infty.$ As $\beta$ changes sign from
positive to negative all the particles start to mutually attract each
others and as $\beta\rightarrow\beta_{N}$ a large number of particles
concentrate along the subvariety of $X$ defined by the center of
the valuation $v_{*}.$

\end{document}